\documentclass[journal]{IEEEtran} 

\usepackage[binary-units=true]{siunitx}
\usepackage{times}
\usepackage{amsmath,amssymb}
\usepackage{accents}
\usepackage{graphicx}
\usepackage{bm}
\let\labelindent\relax
\usepackage[inline]{enumitem}
\usepackage{array}
\usepackage{multicol}
\usepackage{multirow}
\usepackage{tabulary}
\usepackage{booktabs}
\usepackage{subfig}
\usepackage{todonotes}
\usepackage{cite}
\usepackage[hidelinks]{hyperref}
\usepackage{algorithm}
\usepackage{algpseudocode}
\usepackage{breqn}
\usepackage{tcolorbox}
\usepackage{mathtools}
\usepackage{cuted}

\allowdisplaybreaks

\DeclareMathOperator{\diag}{diag}

\DeclareMathOperator*{\mini}{min.}

\definecolor{light-gray}{gray}{0.95}
\definecolor{dark-gray}{gray}{0.5}
\definecolor{mygray}{gray}{0.75}
\newcommand{\BIN}{\begin{bmatrix}}
\newcommand{\BOUT}{\end{bmatrix}}

\pgfdeclarelayer{bg1}   
\pgfdeclarelayer{bg}  
\pgfsetlayers{bg1,bg,main}  

\definecolor{orange}{rgb}{0.99,0.69,0.07}
\definecolor{lightgray}{gray}{0.85}
\definecolor{light-gray}{gray}{0.95}
\definecolor{dark-gray}{gray}{0.5}

\tikzset{cross/.style={cross out, draw=black, minimum size=2*(#1-\pgflinewidth), inner sep=0pt, outer sep=0pt},
cross/.default={1pt}}

 \makeatletter
 \newcommand\fs@spaceruled{\def\@fs@cfont{\bfseries}\let\@fs@capt\floatc@ruled
   \def\@fs@pre{\vspace{5pt}\hrule height.8pt depth0pt \kern2pt}%
   \def\@fs@post{\kern2pt\hrule\relax}%
   \def\@fs@mid{\kern2pt\hrule\kern2pt}%
   \let\@fs@iftopcapt\iftrue}
 \makeatother

\title{$\mathcal{N}$\hspace{-3pt}IPM-MPC: An Efficient Null-Space Method Based Interior-Point Method for Model Predictive Control}
\author{Kai Pfeiffer$^{1}$, Ludovic Righetti$^{1}$%
	\thanks{$^{1}$Machines in Motion Laboratory, Tandon School of Engineering, New York University, New York, USA}%
	\thanks{Part of this work was supported by New York University, NSF grants
		1936332, 1824434, 1833666, 1564142, 1925079, 1825993; NYU WIRELESS and its industrial affiliates; NIST grant 70NANB17H166; SRC; and
		a research grant from OPPO.}
}

\begin{document}
	\maketitle 
	\thispagestyle{empty}
	\pagestyle{empty}
	
	\begin{abstract}
	Linear Model Predictive Control (MPC) is a widely used method to control systems with linear dynamics. Efficient interior-point methods have been proposed which leverage the block diagonal structure of the quadratic program (QP) resulting from the receding horizon control formulation. However, they require two matrix factorizations per interior-point method iteration, one each for the computation of the dual and the primal. Recently though an interior point method based on the null-space method has been proposed which requires only a single decomposition per iteration. While the then used null-space basis leads to dense null-space projections, in this work we propose a sparse null-space basis which preserves the block diagonal structure of the MPC matrices. Since it is based on the inverse of the transfer matrix we introduce the notion of so-called virtual controls which enables just that invertibility. A combination of the reduced number of factorizations and omission of the evaluation of the dual lets our solver outperform others in terms of computational speed by an increasing margin dependent on the number of state and control variables.
	\end{abstract}
	
	\section{Introduction}

	Model Predictive Control aims to efficiently minimize a given cost function with respect to a system's \textit{state} $x$ and \textit{input} or \textit{control} $u$. By virtue of the knowledge of the system's dynamics and therefore its future behavior the cost can be minimized not only instantaneously but such that the current control is also optimal with regard to the system response at some future point. Model Predictive control is therefore also referred to as receding horizon control as the system's behavior at time $t$ is optimized with respect to a certain time window of length $T$ into the future. 
	
	Linear Model Predictive Control (MPC) handles systems with linear dynamics affine both in the state and the control variables. Such problems are common both in industrial control applications such as plants~\cite{Rawlings2009} but also mobile platforms like planes~\cite{Gibbens2011}, cars~\cite{Luo2010} or robots~\cite{Wieber2006}.
	
	In combination with a quadratic objective a MPC optimization problem takes the shape of a quadratic program (QP) with linear equality and inequality constraints. Solving such problems is well established. Active-set methods~\cite{qpOASES} iterate on the so-called active set until the optimal one is found where all active inequality constraints hold as equalities. Active-set methods can be warm-started very efficiently~\cite{Herceg2015} with the previously found active-set and the assumption that it only changes slowly over the control (an assumption which holds well in practice; large changes can be observed for badly posed problems\cite{Pfeiffer2018,pfeiffer2020} however). 
	The Alternating Direction Method of Multipliers (ADMM)~\cite{Gabay1976,Boyd2011} is an operator-splitting method with very good practical convergence behavior, warm-start capabilities and low iteration cost. While the solver is not designed to converge at high accuracy, it has been recently proposed to derive active set guesses from intermediate solver states and solve the corresponding equality only problem~\cite{osqp}. If the active-set guess is correct the solution has zero primal and dual residual. The same framework also overcomes tuning difficulties associated with operator-splitting methods and requires only a single factorization for solving the QP, making it very fast especially for large-scale problems.
	On the contrary, interior-point methods (IPM)~\cite{Nesterov1994} based on the Mehrotra predictor-corrector algorithm~\cite{Mehrotra1992} require little to no tuning effort by design and converge to a high degree. They iterate on a linear approximation of the non-linear optimization problem  and are likewise applicable for problems with large numbers of inequality constraints~\cite{bartlett2000} as it is usually the case for MPC with both the states and the controls being bounded over the horizon. Warm-start strategies for linear programming~\cite{Yildirim2002,Elizabeth2008} and MPC~\cite{Shahzad2010} have been proposed but those capabilities are limited, therefore making it potentially computationally inferior to the ADMM if highly accurate solutions are not required.
	
	Receding horizon control exposes a block diagonal structure where each time step is only coupled with the previous and the next one. IPM solvers are typically capable of exploiting this sparsity~\cite{gurobi,qpSwift2019}. This way the computational complexity of solving MPC's only grows linearly with the length of the horizon and not cubically as it would be the case if a dense QP was solved~\cite{wangboyd2010}. Further extensions have been proposed to compose matrix products more efficiently or update Cholesky factors instead of calculating them anew for low rank inequality constraints~\cite{zeilinger2012}. On the other hand, the authors in~\cite{Frison2016} propose a solver which works on a condensed version of the sparse MPC problem and therefore has cubic time complexity in the horizon length. This can be advantageous since linear algebra routines work more efficiently on large dimensions.
	
	The two sparse MPC solvers based on the Schur complement method~\cite{wangboyd2010,zeilinger2012} have one commonality: two decompositions are necessary to solve a single iteration of the IPM, one to obtain the dual associated with the dynamics constraints and one to obtain the primal. Recently though the solver NIPM-HLSP based on the null-space method has been proposed which only requires a single decomposition per iteration~\cite{pfeiffer2021}. Since the obtained primal step is in the null-space of the dynamic constraints their feasibility is always ensured (given that the initial primal is feasible), rendering the calculation of the dual unnecessary and in the consequence reducing the operational count additionally.
	
	However, NIPM-HLSP was designed for instantaneous control without the characteristic block diagonal structure of receding horizon control problems. Consequently, the use of a dense basis of the null-space was perfectly valid. However, in the MPC case such a basis would destroy the sparsity and the computational effort of resolving the MPC would now increase cubically with the length of the receding horizon which is prohibitive. While the null-space method for QP problems in MPC has been applied for example in~\cite{HPIPM}, the authors do not further specify a sparsity maintaining null-space basis.
	
	Our contributions are therefore two-fold:
	\begin{itemize}
		\item We propose a sparse null-space basis which preserves the block diagonal sparsity of the matrix to decompose. 
		\item This null-space basis assumes the invertibility of the transfer matrix. Since this is usually not the case we introduce the concept of so-called virtual control which achieves just that invertibility.
	\end{itemize}

	Our proposed solver based on these two contributions outperforms other available linear MPC solvers based on the IPM by up to 70\% per Newton iteration depending on the problem formulation.
	
	The paper is composed as follows: in section~\ref{sec:lmpcqp} we recall the formulation of receding horizon problems as QP's with special block-diagonal structure. In section~\ref{sec:nsipm} we recall the null-space method based IPM which requires only a single decomposition per solver iteration instead of two. In section~\ref{sec:nsbasis} a sparse null-space basis is proposed which maintains the MPC sparsity of the projected matrices. This null-space basis requires the invertibility of the transfer matrix which we ensure by introducing so-called virtual controls in the case of under-actuated systems, see sec.~\ref{sec:virtualCtrls}. Section~\ref{sec:sparseQR} presents sparse QR decompositions for the special case of a square transfer matrix. In sec.~\ref{sec:opComp} we compare the operational counts between the different IPM solvers. The evaluation (sec.~\ref{sec:evaluation}) confirms that our solver formulation indeed can be faster than IPM formulations based on the classical normal equations. Finally, we conclude the paper in sec.~\ref{sec:conclusion} with some thoughts on potential future work.
	
	\section{Linear Model Predictive control as a QP}
	\label{sec:lmpcqp}
	
	Linear model predictive control is concerned with systems with linear dynamics of the discretized form 
	\begin{equation}
	x(t + 1) = A_{x,e}x(t) + B_{u,e}u(t) + c(t) \qquad t = 0,1,\dots
	\end{equation}
	$t$ denotes time, $x\in\mathbb{R}^{n_x}$ is the state vector, $u\in\mathbb{R}^{n_u}$ is the control input and $c$ is a disturbance.
	$A_{x,e}\in\mathbb{R}^{n_x,n_x}$ and $B_{u,e}\in\mathbb{R}^{n_x,n_u}$  are the system and transfer matrices, respectively. Since the system and transfer matrix do not change over time such dynamics are referred to as \textit{time-invariant}. Throughout this paper we assume both these matrices to be full rank, i.e. there are no linear dependent matrix rows or columns.
	
	The optimization problem associated with a receding horizon control problem now takes following form
	\begin{align}
	\mini_{x,u} \qquad & J\coloneqq l_f(x(t+T)) + \sum_{\tau = t}^{t + T - 1}l(x(\tau),u(\tau))\\
	\text{s.t}
	\qquad&	x(\tau + 1) = A_{x,e}x(\tau) + B_{u,e}u(\tau) + w\nonumber\\
	& A_{x,i}x(\tau) + B_{u,i}u(\tau) \geq b_{x,u,i}\nonumber\\
	&\tau = t, \dots, t + T - 1\nonumber
	\end{align}
	Goal of this problem is to minimize the cost function $l(x(t),u(t))$ over the receding horizon of length $T$ while respecting the dynamics and given inequality constraints $A_{x,i}\in\mathbb{R}^{m_i,n_x}$ and $B_{u,i}\in\mathbb{R}^{m_i,n_u}$. We specifically consider quadratic cost functions of the form
	\begin{equation}
	l(x,u) = \BIN
	x\\u
	\BOUT^T
	\BIN
	Q & S\\
	S^T & V
	\BOUT
	\BIN
	x\\
	u
	\BOUT
	+
	q^Tx + r^Tu
	\end{equation}
	The weight matrix and its constituting components $Q\in\mathbb{R}^{n_x,n_x}$ and $U\in\mathbb{R}^{n_u,n_u}$ are positive definite. $S\in\mathbb{R}^{n_x,n_u}$ is some cross term coupling the cost of the states and controls.
	We omitted the time index $t$ for better readability and do so throughout this paper if applicable.

	The above problem is a QP. We rewrite it to
	\begin{align}
	\mini_{y} \qquad & y^THy + g^T y
	\label{eq:QP}\\
	\text{s.t}
	\qquad& A_e y - b_e	= 0	\nonumber\\
	\qquad& A_{i} y - b_{i}\geq 0\nonumber
	\end{align}
	with $y=\BIN u(t)^T,x(t + 1)^T,\dots,u(T-1),x(T)^T\BOUT^T\in\mathbb{R}^{n}$ ($n=T(n_x+n_u)$),
	$H\in\mathbb{R}^{n,n}$, $g\in\mathbb{R}^{n}$, $A_e\in\mathbb{R}^{Tn_x,n}$, $b_e\in\mathbb{R}^{Tn_x,n}$, $A_i\in\mathbb{R}^{Tm_i,n}$ and $b_i\in\mathbb{R}^{Tm_i,n}$ defined as follows 
	\begin{align}
	&H \coloneqq 
	\BIN 
	U & 0 & 0 & \cdots & 0 & 0 & 0 \\
	0 & Q & S & \cdots & 0 & 0 & 0\\
	0 & S^T & U & \cdots & 0 & 0 & 0 \\
	\vdots & \vdots & \vdots & \ddots & \vdots & \vdots & \vdots \\
	0 & 0 & 0 & \cdots & Q & S & 0\\
	0 & 0 & 0 & \cdots & S^T & U & 0\\
	0 & 0 & 0 & \cdots & 0 & 0 & Q_f
	\BOUT\\
	&g \coloneqq 
	\BIN (r + 2S^Tx(t))^T&
	q^T&
	r^T&
	\cdots&
	q^T&
	r^T&
	q_f^T
	\BOUT^T\\
	&A_e \coloneqq\\ 
	&
	\resizebox{0.475\textwidth}{!}{$
	\BIN 
	-B_{u,e} & I & 0 & 0 & \cdots & 0 & 0 & 0\\
	0 & -A_{x,e} & -B_{u,e} & I & \cdots & 0 & 0 & 0\\
	0 & 0 & 0 & -A_{x,e} & \cdots & 0 & 0 & 0\\
	\vdots & \vdots & \vdots & \ddots & \vdots & \vdots & \vdots \\
	0 & 0 & 0 & 0 & \cdots & I & 0 & 0\\
	0 & 0 & 0 & 0 & \cdots & -A_{x,e} & -B_{u,e} & I\\
	\BOUT
	$}\nonumber\\
	&b_e \coloneqq
	\BIN
	(A_{x,e}x(t) + c(t))^T&
	c(t + 1)^T&
	\cdots&
	c(T)^T
	\BOUT^T\\
	&A_i \coloneqq 
	\BIN B_{u,i} & 0 & 0 & \cdots & 0 & 0 & 0\\
	0 & A_{x,i} & B_{u,i} & \cdots & 0 & 0 & 0\\
	\vdots & \vdots & \vdots & \ddots & \vdots & \vdots & \vdots \\
	0 & 0 & 0 & \cdots & A_{x,i} & B_{u,i} & 0\\
	0 & 0 & 0 & \cdots & 0 & 0 & A_{x,i}
	\BOUT\\
	&b_i = 
	\BIN
	(b_{xu,i} - A_{x,i}x(t))^T &
	b_{xu,i}^T&
	\cdots&
	b_{xu,i}^T&
	b_{xu,i,f}^T
	\BOUT^T
	\end{align}

	\section{The null-space method based primal-dual IPM}
	\label{sec:nsipm}
	
	The above QP~\eqref{eq:QP} can solved by the IPM.	
	We apply the IPM by introducing the slack variable $w_i$ and bounding it away from zero by the log-barrier function
	\begin{align}
	\mini_{x} \qquad & \frac{1}{2}y^THy + g^T y - \sigma\mu\sum \log(w_i)\\
	\text{s.t}
	\qquad& A_e y - b_e	= 0	\nonumber\\
	\qquad& A_{i} y - b_{i} = w_i\nonumber\\
	\qquad& c_i \geq 0\nonumber
	\end{align}
	
	The Lagrangian of this optimization problem writes as
	\begin{equation}
	\begin{aligned}
	\mathcal{L} &= \frac{1}{2}y^THy + g^Ty -\sigma\mu\sum\log(w_{i}) \\ 
	&- \lambda_{e}^T(A_{e}y - b_{e})
	- \lambda_{i}^{T}(A_{i}y - b_{i} - w_{i}) 
	\end{aligned}
	\end{equation}
	$\lambda_{e}$ are the Lagrange multipliers associated with the equality constraints and $\lambda_{i}$ are the ones associated with the inequality constraints.
	
	The duality measures $\mu$ and the centering parameters $\sigma$~\cite{vanderbei2013} are given by
	\begin{align}
	\mu &= \lambda_{i}^Tw_{i}/(n_z+{m}_{i}) \qquad \text{and} \qquad
	\sigma\in[0,1]
	\end{align}
	or the values can be determined by Mehtrotra's predictor-corrector algorithm~\cite{Mehrotra1992}.
	
	The minimizer is found at the Lagrangian's stationary points $\nabla_q \mathcal{L} = 0$. 		$q$ is the variable vector
	\begin{equation}
	q \coloneqq \BIN y^T & \lambda_{e}^T & \lambda_{i}^T  & w_{i}^T \BOUT^T
	\end{equation}
	We get the slightly rewritten KKT conditions
	\begin{align}
	K
	=
	\BIN
	k_1\\
	k_2\\
	k_3\\
	k_4
	\BOUT&\coloneqq
	\BIN
	Hy + g -A_{e}^T\lambda_{e} -A_{i}^T\lambda_{i} \\
	b_{e} - A_{e}y\\
	b_{i} - A_{i}y + w_{i}\\
	\lambda_{i}\odot w_{i} - \sigma\mu e
	\BOUT
	=
	0\label{eq:k}
	\end{align}
	$\odot$ is the element-wise product of two vectors.
	$W_{i} = \diag(w_{i})$ and $\Lambda_{i} = \diag(\lambda_{i})$ are square matrices with the vectors $w_{i}$ and $\lambda_{i}$ as diagonals. $e\in\mathbb{R}^{m_{i}}$ is a vector of ones.
	
	We additionally have the feasibility conditions
	\begin{equation}
	w_{i} \geq 0 \qquad \text{ and } \qquad \lambda_{i} \geq 0
	\label{eq:feasCond}
	\end{equation}
	
	Finally, we linearize this nonlinear equation by the Newton step
	\begin{align}
	K(q + \Delta q) = K(q) + \nabla_q K \Delta q=0
	\end{align}
	with
	\begin{equation}
	\nabla_qK=
	\BIN
	H & -A_{e}^T & -A_{i}^T &  0\\
	-A_{e} & 0 & 0 & 0 \\
	-A_{i} & 0 & 0 & I \\
	0 & 0 & W_{i} & \Lambda_{i}
	\BOUT
	\end{equation}

	This linear equation is now solved iteratively and a new step $\alpha\Delta q$ is applied to $q$ in each iteration. $\alpha$ is determined by line search and ensures the feasibility conditions~\eqref{eq:feasCond}.
	
	For a more efficient algorithm we can apply substitutions for $\Delta w_{1,i}$  and $\Delta \lambda_{i}$
	which yields the augmented system~\cite{zeilinger2012}
	\begin{align}
	\BIN
	\Phi & -A_{e}^T\\
	-A_{e} & 0
	\BOUT
	\BIN
	\Delta x\texttt{}\\
	\Delta \lambda_{e}
	\BOUT
	=
	\BIN
	r_1\\
	r_2
	\BOUT
	\label{eq:qpNeNmethodFull}
	\end{align}
	The matrix on the left is square, symmetric and indefinite. 
	\begin{equation}
	\Phi = H + A_{i}^TW_{i}^{-1}\Lambda_{i} A_{i}
	\end{equation}
	is positive definite.
	The right hand side is given by
	\begin{align}
	r_1 &= -k_1 + A_{i}^TW_{i}^{-1}(\Lambda_{i}k_3 - k_{4}) \\
	&=-Hx - g + A_{e}^T\lambda_{e} +A_{i}^TF\nonumber\\
	r_2 &= -k_{2} =  A_{e}x - b_{e}
	\end{align}
	with
	\begin{align}
	F = \lambda_{i}
	+ W_{i}^{-1}(\lambda_{i} \cdot (b_{i} - A_{i}x))
	\label{eq:fpredictor}
	\end{align}
	for the predictor step
	and 
	\begin{align}
	F = \lambda_{i}
	+ W_{i}^{-1}(\lambda_{i} \cdot (b_{i} - A_{i}x) - \Delta\lambda_i\cdot\Delta w_i+ \sigma\mu e)
	\label{eq:fcorrector}
	\end{align}
	for the corrector step of Mehrotra's predictor-corrector algorithm.
	
	Solving the augmented system directly is not advised since this would ignore the given sparsity in the bottom right corner. Instead, as proposed in~\cite{zeilinger2012}, the Schur complement can be formed such that we obtain the \textit{classical normal equations}. First, the dual is computed by
	\begin{align}
		A_e\Phi^{-1}A_e^T\Delta \nu &= -r_2 - A_eC^{-1}r_1
		\label{eq:clasNEq}
	\end{align}
	and with it the primal is recovered
	\begin{align}
		\Delta x &= \Phi^{-1}(r_1 + A_e^T\Delta\nu)
	\end{align}
	As can be seen, two Cholesky decompositions have to be computed per Newton iteration, one of $\Phi$ and one of $A_e\Phi^{-1}A_e^T$. As recently proposed~\cite{pfeiffer2021} this can be prevented by using the null-space method. 
	First, we assume that 
	\begin{equation}
	r_2 = 0
	\label{eq:feasPt}
	\end{equation}
	by obtaining an initial feasible $x$ by solving $A_{e}x = b_{e}$ beforehand, for example with a sparse QR decomposition of $A_e$ as described in sec.~\ref{sec:sparseQR}. We then apply the null-space method~\cite{Nocedal2006} by introducing the variable change
	\begin{equation}
	\Delta x = N\Delta z
	\label{eq:changeOfVarQP}
	\end{equation}
	$N$ is a basis of the null-space of $A_{e}$ such that $A_{e}N = 0$. Consequently, $A_{e}\Delta x = 0$ such that the condition $r_2 = 0$ continues to be fulfilled.
	By furthermore projecting the augmented system into the null-space basis $N$ of the equality constraints $A_e$ we get the \textit{projected normal equations}
	\begin{equation}
	N^T\Phi N\Delta z =N^Tr_1 
	\label{eq:projNEq}
	\end{equation}
	Now only one decomposition of a reduced system of less variables needs to be conducted per Newton iteration. This is due to the concept of so-called variable elimination~\cite{dimitrov:2015} reducing the number of variables $n$ by the rank of the equality constraints $Tn_x$ such that $N^T\Phi N\in\mathbb{R}^{n - Tn_x = Tn_u}$ with $N\in\mathbb{R}^{2Tn_x,Tn_x}$. With our choice of variable-reducing null-space basis, which maintains the block-diagonal structure of $\Phi$ (see sec.~\ref{sec:nsbasis}), the full-rank and positive definite projected normal equations can be solved by a block-wise Cholesky decomposition.
	
	Note that if the quadratic objective~\eqref{eq:QP} is of least-squares form the projected normal equations can also be expressed in least-squares form~\cite{pfeiffer2021}. This can be advantageous if the number of inequality constraints is small. Since this is not the case for most control applications we will not further address it.
	
	The dual step $\Delta\lambda_e$ is computed by solving
	\begin{align}  
	A_e^T \Delta\lambda_e = \Phi\Delta x - r_1
	\label{eq:dualStep}
	\end{align}
	
	As observed in~\cite{pfeiffer2021}, the Lagrange multipliers $\Delta\lambda_e$ do not necessarily need to be computed  since none of the other other primal and dual variables depend on it. They are only necessary for the evaluation for the KKT vector and the convergence criteria $\left\Vert K \right\Vert < \epsilon$ with the small numerical threshold $\epsilon = 10^{-9}$. This is in contrast to the classical normal equations which require the evaluation of the dual in each Newton iteration in order to obtain the primal. 
	We base our convergence test additionally on the norm of the residual $\left\Vert N^Tr_1 \right\Vert^2< \epsilon$, making the evaluation of the Lagrange multipliers obsolete after all. This is valid since the null-space method maintains primal feasiblity with respect to the equality constraints and we do not further need to evaluate $k_1$ and $k_2$.

	\section{A sparse null-space basis}
	\label{sec:nsbasis}
	
	In order to apply the null-space method we require a basis of the null-space of $A_e$. In~\cite{pfeiffer2021} a basis based on the QR decomposition of $A_e$ was used. It is straightforward to compute and can handle rank deficiencies in $A_e$. However, the basis is dense and therefore any projection is dense, too. Especially in the context of MPC this would be highly disadvantageous since then the computational effort of decomposing $N^T\Phi N$ would grow cubically with the length of the receding horizon $T$. This would be in contrast with the approaches in~\cite{wangboyd2010,zeilinger2012} which preserve sparsity and whose operation counts only grow linearly with the horizon length.
	
	It is therefore important to identify a sparse null-space basis preserving the block diagonal structure of the projection. 
	We assume that both $A_{x,e}$ and $B_{u,e}$ are full column-rank (in the sense that there are no linear dependent rows or columns). 
	Such a sparse basis of the null-space is given by
	\begin{align}
	N = \BIN
	I & 0 & \cdots & 0 \\
	B_{u,e} & 0 & \cdots & 0\\
	C & I & \cdots & 0\\
	0 & B_{u,e} & \cdots & 0\\
	0 & C & \cdots & 0\\
	\vdots & \vdots & \ddots & \vdots\\
	0 & 0 & \cdots & I\\		
	0 & 0 & \cdots & B_{u,e}
	\BOUT
	\in\mathbb{R}^{2Tn_x,Tn_u}
	\end{align}
	with $C = -B_{u,e}^{-1}A_{x,e}B_{u,e}$. 	
	We see that if $B_{u,e}\in\mathbb{R}^{n_x,n_u}$ is invertible and $B_{u,e}^{-1}$ exists (in case of a full-rank square matrix; right inverse $B_{u,e}^{+}$ for $n_u>n_x$) the projection of each instance of the horizon is only coupled with the one from the previous and the next one. The invertibility can be achieved for the common control constellation $n_u < n_x$ by virtue of a concept we refer to as \textit{virtual controls} and which is further detailed in the next section~\ref{sec:virtualCtrls}. 
	
	The single entries of the projection $A_eN$ take the form of either $-B_{u,e} + B_{u,e} = 0$, $-A_{x,,e}B_{u,e} + B_{u,e}B_{u,e}^{-1}A_{x,e}B_{u,e}=0$ or simply zero due to the sparsity of the null-space basis. The projections $N^THN$ and $N^TA_i^TA_iN$ are given in the following:
	\begin{align}
	&N^THN =\\
	& \BIN
	R + B_{u,e}^TM_1 + C^TM_2& M_2^T & \cdots & 0\\
	M_2 & R + \cdots& \cdots & 0\\
	0 & M_2 & \cdots & 0\\
	\vdots & \vdots & \ddots & \vdots\\
	0 & 0 & \cdots & M_2^T\\
	0 & 0 & \cdots & R + B_{u,e}^TQ_fB_{u,e}
	\BOUT
	\label{eq:nhn}
	\end{align}
	with 
	\begin{align}
		M_1 = QB_{u,e} + SC \qquad \text{and} \qquad
		M_2 = S^TB_{u,e} + RC
	\end{align}	
	and

	\vspace{20pt}
	\begin{strip}
	\begin{align}
	&N^TA_i^T\Xi_iA_iN = \label{eq:nain}\\
	&
		\BIN
	B_{u,i}^T\Xi_i(t)B_{u,i}+ M_3^T\Xi(t + 1)M_3 & M_3^T\Xi_i(t + 1)B_{u,i} & \cdots & 0 \\
	B_{u,i}^T\Xi_i(t + 1)M_3 & B_{u,i}^T\Xi_i(t + 1)B_{u,i} + \cdots & \cdots & 0\\
	\vdots & \vdots & \ddots & \vdots \\
	0 & 0 & \cdots & M_3^T\Xi_i(T - 1)B_{u,i}\\
	0 & 0 & \cdots & B_{u,i}^T\Xi_i(T-1)B_{u,i} + B_{u,e}^TA_{x,i}^T\Xi_i(T)A_{x,i}B_{u,e}
	\BOUT
		\nonumber
		\end{align}
		\end{strip}

	with $\Xi_i(t) = W^{-1}(t)\Lambda(t)$ and $M_3 = A_{x,i}B_{u,e} + B_{u,i}C.$
	
	As can be easily seen, the chosen sparse null-space basis maintains the block diagonal sparsity of the MPC matrices which takes the symbolic form
	\begin{align}
		N^T\Phi N &= N^T (H + A_i^T\Xi_iA_i)N\\
		&=\BIN
		Y_{11} & Y_{12} & 0 & \cdots & 0 & 0\\
		Y_{21} & Y_{22} & Y_{23} & \cdots & 0 & 0 \\
		0 & Y_{32} & Y_{33} & \cdots & 0 & 0 \\
		\vdots & \vdots & \vdots & \ddots & \vdots & \vdots\\
		0 & 0 & 0 & \cdots & Y_{T-1,T-1} & Y_{T-1,T}\\
		0 & 0 & 0 & \cdots & Y_{T,T-1} & Y_{TT}
		\BOUT\nonumber
	\end{align}
	where the single entries $Y_{ij}$ follow from the addition of~\eqref{eq:nhn} and~\eqref{eq:nain}. The block-wise Cholesky decomposition $Y = LL^T$ with the lower bidiagonal block matrix
	\begin{equation}
		L = \BIN
		L_{11} & 0 & 0 & \cdots & 0 & 0 \\
		L_{21} & L_{22} & 0 & \cdots & 0 & 0\\
		0 & L_{32} & L_{33} & \cdots & 0 & 0 \\
		\vdots & \vdots & \vdots & \ddots & \vdots & \vdots\\
		0 & 0 & 0 & \cdots & L_{T - 1,T-1} & 0 \\
		0 & 0 & 0 & \cdots & L_{T,T-1} & L_{TT}
		\BOUT
	\end{equation}
	is then computed as described in~\cite{wangboyd2010}
	\begin{align}
		&L_{11} \leftarrow \text{CHOL}(Y_{11})\\
		&\text{for }i=2,\dots,T:\nonumber\\
		&\qquad L_{i,i-1}\leftarrow L_{i,i-1}L_{i-1,i-1}^T = Y_{i,i-1} \\
		&\qquad L_{ii} \leftarrow \text{CHOL}(Y_{ii} - L_{i,i-1}L_{i,i-1}^T)
	\end{align}
	$\text{CHOL}$ is an in-place Cholesky decomposition routine and $L_{i,i-1}$ is obtained by solving the corresponding equation by backward substitution.

	\section{Notion of virtual controls}
	\label{sec:virtualCtrls}
	
	In a reasonable control context we usually have $n_u < n_x$ (i.e. underactuated systems), meaning that $B_{u,e}$ is not invertible. In this case we can only calculate a left pseudo-inverse which does not fulfill the necessary condition $B_{u,e}B_{u,e}^+ = I$. However, we can achieve invertibility of the transfer matrix by introducing so-called \textit{virtual controls} $u^*\in\mathbb{R}^{n_{u^*}}$ of dimension $n_{u^*} = n_x - n_u$. This lets us obtain the new full-rank square and therefore invertible transfer matrix
	\begin{equation}
	B_{\hat{u},e} = \BIN B_{u,e} & B_{u^*,e} \BOUT \in\mathbb{R}^{n_x,n_x}
	\end{equation}
	$\hat{u}$ is the new control vector consisting of the original and virtual controls $u$ and $u^*$ as $\hat{u}\coloneqq\BIN u^T & u^{*T} \BOUT^T$ with $n_{\hat{u}} = n_x$.
	$B_{u^*,e}$ is defined as
	\begin{equation}
	B_{u^*,e} = Q_2\diag(r_{min})
	\label{eq:rmin}
	\end{equation}
	where $Q_2$ is taken from the QR decomposition of $B_{u,e}$
	\begin{equation}
	B_{u,e} = \BIN Q_1 & Q_2 \BOUT \BIN R \\ 0 \BOUT
	\end{equation}
	$r_{min}$ is the absolute smallest entry on the diagonal of $R$ in order to maintain a well conditioned matrix $B_{\hat{u},e}$. 
	
	The QR decomposition of $B_{\hat{u},e}$ is then
	\begin{align}
	B_{\hat{u},e} &= Q_{\hat{u},e}R_{\hat{u},e} = \BIN Q_1 & Q_2 \BOUT \BIN R & 0 \\ 0 & \diag(r_{min}) \BOUT 
	\label{eq:bhatueqr}\\
	&\left(= \BIN Q_1 R & Q_2\diag(r_{min}) \BOUT = \BIN B_{u,e} & B_{u^*,e} \BOUT \right)\nonumber
	\end{align}
	$C$ can be calculated as
	\begin{equation}
		C = -R_{\hat{u},e}^{-1}Q_{\hat{u},e}^{T}A_{x,e}B_{\hat{u},e}
	\end{equation}
	The new variable vector of the QP is $y=\left[ u(t)^T,u^*(t)^{T},x(t \hspace{-2pt}+\hspace{-2pt} 1)^T,\dots,u(T\hspace{-2pt}-\hspace{-2pt}1)^T,u^*(T\hspace{-2pt}-\hspace{-2pt}1)^{T},x(T)^T\right]^T \\\in \mathbb{R}^{n}$ with $n=2Tn_x$ entries. 
	The projected matrix is now of dimension $N^T\Phi N\in\mathbb{R}^{Tn_x,Tn_x}$ after we added $Tn_{u^*}$ but eliminated $Tn_x$ variables. 
	Consequently, the projected matrix $N^T\Phi N$ is smaller or equal in dimension with respect to the unprojected matrix $\Phi\in\mathbb{R}^{T(n_x + n_u),T(n_x + n_u)}$ depending on the input size $n_u$ (equal for for the trivial case $n_u = 0$, smaller for $n_u>0$). 
	A detailed computational comparison is given in sec.~\ref{sec:opComp}.
	
	For good measure, we choose the new control cost weight matrix to be
	\begin{equation}
		R\coloneqq\BIN R & 0 \\ 0 & I \BOUT\in\mathbb{R}^{n_{\hat{u}},n_{\hat{u}}}
	\end{equation}
	or some well conditioned variant of it with a weighted identity matrix. Since the projection $N^THN$~\eqref{eq:nhn} is full-rank because of the term $B_{\hat{u},e}QB_{\hat{u},e}$ and therefore applicable for the Cholesky decomposition regardless this is not necessarily required.
	
	The remaining equations and definitions from the previous sections hold as long as the original control vector $u$ is replaced by the new one $\hat{u}$.
	
	The virtual controls are not allowed to change the behavior of the original dynamic system so we have to introduce the equality constraints
	\begin{equation}
	u^* = 0
	\end{equation}
	in order to eliminate any effects of the virtual control.
	However, including these into the set of equality constraints renders $B_{\hat{u},e}$ again into a non-invertible rectangular matrix $\in\mathbb{R}^{n_x + n_{u^*},n_x}$. In order to circumvent this we reformulate our virtual control constraints into a set of inequality constraints
	\begin{align}
	u^* \leq 0 \qquad\text{and}\qquad
	u^* \geq 0
	\end{align}
	and add them to the set of inequality constraints $B_{\hat{u},i}$. This may raise questions about the numerical stability of the IPM but in practice we did not observe any such adverse effects. This can be attributed to our primal-dual IPM formulation which only maintains dual $w_i\geq 0$ but not necessarily primal feasibility $A_iy - b_i\ngeq 0$ throughout the Newton iterations. In the same line of argument, the concept of virtual controls may not be realizable in primal-barrier formulations of the IPM~\cite{wangboyd2010}.
	
	\section{A sparse QR decomposition of $A_e$}
	\label{sec:sparseQR}
	
	In order to obtain the initial feasible equality point $A_ex=b_e$~\eqref{eq:feasPt} we require the QR decomposition of $A_e$. The QR decomposition of $A_e$ has to be permuting. In order to determine the permutation order we first identify the approximate condition number $\tilde{\kappa}$ of $A_{x,e}$ and $B_{\hat{u},e}$ by calculating the ratio of the largest to smallest value on the diagonal of their upper triangular QR factor $R_{x,e}$ and $R_{\hat{u},e}$. We choose its threshold $\xi = 10$. Note that the amplification factor during the block-wise inversion is of the order $\xi^T$ and may require further tuning depending on the problem at hand. Since this operation is done offline we could also obtain the exact condition number but found this criteria to be sufficient in practice. With the assumption of a square and full-rank transfer matrix $B_{\hat{u},e}$ (as is it in our case, see sec.~\ref{sec:virtualCtrls}) we choose following permutation orders and their associated QR decompositions:
	\begin{itemize}
		\item If $\tilde{\kappa}(B_{\hat{u},e})<\xi$ 
		\begin{align}
	&Q\BIN R & V \BOUT \\
	=& A_eP_1 = \BIN \diag\left(\BIN -B_{\hat{u},e} & \dots & -B_{\hat{u},e}\BOUT\right) & D \BOUT\nonumber
	\end{align}
	with the QR decomposition
	\begin{align}
	Q &= \diag\left(\BIN Q_{\hat{u},e} & \dots & Q_{\hat{u},e} \BOUT\right) \\
	R &= \diag\left(\BIN R_{\hat{u},e} & \dots & R_{\hat{u},e} \BOUT\right)
	\end{align}
	This decomposition is the cheapest since it only requires the QR decomposition of $B_{\hat{u},e}$~\eqref{eq:bhatueqr} and no off-diagonal elements need to be handled during its block inversion of cost $O(T2n_x^2)$.
	\item If $\tilde{\kappa}(B_{\hat{u},e})\geq\xi$ but $\tilde{\kappa}(A_{x,e})<\xi$
	\begin{align}
	&Q\BIN R & V \BOUT = A_eP_2 \\
	&=
	\BIN -B_{\hat{u},e} & I & 0 & \cdots & 0 & \\
	0 & -A_{x,e} & I &  \cdots & 0 & \\
	0 & 0 & -A_{x,e} & \cdots & 0 & D\\
	\vdots & \vdots & \vdots & \ddots & \vdots  &  \\
	0 & 0 & 0 & \cdots & I &  & \\
	0 & 0 & 0 & \cdots & -A_{x,e} &
	\BOUT\nonumber
	\end{align}
	with the QR decomposition
		\begin{align}
	Q &= 
	\BIN
	Q_{\hat{u},e} & 0 & \cdots & 0\\
	0 & Q_{x,e} & \cdots & 0\\
	\vdots & \vdots & \ddots & \vdots \\
	0 & 0 & \cdots & Q_{x,e}
	\BOUT\\
	R &= 
	\BIN
	R_{\hat{u},e} & Q_{u,e}^T & 0 & \cdots & 0\\
	0 & R_{x,e} & Q_{x,e}^T & \cdots & 0\\
	\vdots & \vdots & \vdots & \ddots & \vdots \\
	0 & 0 & 0 & \cdots & Q_{x,e}^T\\
	0 & 0 & 0 & \cdots & R_{x,e}
	\BOUT
	\end{align}
	This decomposition requires QR decompositions of $A_{x,e} = \BIN Q_{x,e} & R_{x,e}\BOUT$ and $B_{\hat{u},e}$~\eqref{eq:bhatueqr}. During the block inversion off-diagonal elements need to be handled, making it slightly more expensive at $O(T3n_x^2)$.
	\item Otherwise
		\begin{align}
	&Q\BIN R & V \BOUT = A_eP_3 \\
	=&
	\BIN I & 0 & \cdots & 0 & 0 &\\
 -A_{x,e} & I &  \cdots & 0 & 0 &\\
 0 & -A_{x,e} & \cdots & 0 & 0 & D\\
 \vdots & \vdots & \ddots & \vdots  &  \vdots &\\
 0 & 0 & \cdots & I & 0 \\
 0 & 0 & \cdots & -A_{x,e} & I\\
\BOUT\nonumber
	\end{align}
	This QR decomposition can not be composed block-wise but needs to be calculated from scratch requiring $O(T8n_x^34/3)$ operations. The block-inversion comes at a cost of $O(T6n_x^2)$.
\end{itemize}
The permutation of $D$ does not need to be further specified. In the same sense, $V = Q^TD$ does not need to be computed since it is not used. This might be advantageous in case of time-variant systems where the sparse QR decomposition of $A_e$ needs to be recomputed in each control iteration.

	The QR decomposition can be reused to calculate the dual step~\eqref{eq:dualStep} if required since it does not matter whether the QR decomposition of $A_e$ or $A_e^T$ is used~\cite{dimitrov:2015}. Conveniently, in the context of MPC with its block diagonal structure it is computationally more efficient to use the QR decomposition of $A_e$ since both the resulting triangular factor $R$ and the Householder matrix $Q$ are sparse. This is in contrast to the QR decomposition of $A_e^T$ which results in a dense $Q$ (there is only an upper Hessenberg triangular part to the left).
	
	\section{Operational comparison}
	\label{sec:opComp}
	
	\begin{figure*}[t!]
		\begin{tcolorbox}
			\begin{minipage}[t]{0.46\textwidth}
				\textbf{Projected normal equations (NIPM-MPC)~\eqref{eq:projNEq}}\\
				\textbf{Offline:}
				\begin{enumerate}
					\item Calculate the sparse QR decomposition of $A_e$ ($P_1$: $O(n_x^34/3)$, $P_2$: $O(n_x^34/3)$, $P_3$: $O(T8n_x^34/3)$)
					\item Do the projections $A_iN$ $O(Tm_i(2n_x)^2)$ and $N^THN$ $O(6Tn_x^3)$
				\end{enumerate}
				\textbf{Online, once per Newton's method:}
				\begin{enumerate}
					\item Calculate the equality point $A_ex = b_e$ 
					\item Calculate the Lagrange multipliers $\lambda_e$\\
					$2\times$ $P_1$: $O(T2n_x^2)$, $P_2$: $O(T3n_x^2)$, $P_3$: $O(T6n_x^2)$
				\end{enumerate}
				\textbf{Online, once per Newton iteration:}
				\begin{enumerate}
					\item Form the matrix product $(A_{i}N)^TW_{i}^{-1}\Lambda_{i}(A_{i}N)$, $O(Tm_{i}n_x^2)$  
					\item Cholesky decomposition of $N^T\Phi N$, $O(Tn_x^34/3)$. 
					\item Solve in $O(T2n_x^2)$ for $\Delta z$ and project for $\Delta x$ in $O(T2n_x^2)$ operations (2$\times$ for Mehrotra's predictor-corrector algorithm)
				\end{enumerate}
				$\sum:$ $O(T((8 + m_{i})n_x^2 + n_x^34/3))$
			\end{minipage}
					\hspace{0.5cm}
			\begin{minipage}[t]{0.46\textwidth}
				\textbf{Classical normal equations (FORCES, case B~\cite{zeilinger2012})~\eqref{eq:clasNEq}}\\
				\textbf{Online, once per Newton iteration:}
				\begin{enumerate}
					{\color{gray}\item Form the matrix product $A_{i}^TW_{i}^{-1}\Lambda_{i}A_{i}$, $O(Tm_{i}(n_x+n_u)^2)$} 
					{\color{gray}\item Cholesky decomposition of $\Phi$,  $O(T(n_x + n_u)^3/3)$}
					\item Form $A_{e}\Phi^{-1}A_{e}^T$ in $O(T4(n_x + n_u)n_x^2)$
					\item Block-wise Cholesky decomposition of $A_{e}\Phi^{-1}A_{e}$, $O(Tn_x^34/3)$
					\item Calculate the dual $\lambda_{e}$ in $O(T2n_x^2)$ operations
					and the primal in $\Delta y$, $O(T2(n_x + n_u)^2)$ operations (2$\times$ for Mehrotra's predictor-corrector algorithm)
				\end{enumerate}
				$\sum:$ $O(T(4(n_x + n_u)n_x^2  + 4(n_x^2 + (n_x+n_u)^2) + n_x^34/3))$\\
			\end{minipage}%
		\end{tcolorbox}
		\caption{Number of operations for the classical and the projected normal equations.}
		\label{fig:NewtonIterOp}
	\end{figure*}

	In fig.~\ref{fig:NewtonIterOp} we give a brief overview of the steps and the number of operations associated both with the classical and projected normal equations. We compare our solver NIPM-MPC to FORCES~\cite{zeilinger2012} which is based on the classical normal equations. We do this for the case of a diagonal matrix $H$ which corresponds to the B variant of the FORCES solver.
	
	Both the QR decomposition of $A_e$ and the projections into the null-space can be done offline in the case of our solver. The calculation of the initial feasible equality point requiring matrix vector multiplications only needs to be done once before the Newton's method. Each iteration then first requires forming the matrix product $(A_{i}N)^TW_{i}^{-1}\Lambda_{i}(A_{i}N)$ with the updated dual $W_{i}^{-1}\Lambda_{i}$. The resulting matrix is dense even if $A_i$ only consists of bound constraints. This is one of the disadvantages of the null-space method based IPM~\cite{pfeiffer2021}. Furthermore required are one Cholesky decomposition of $N^T\Phi N$ and the subsequent solution of the linear KKT system. The primal step $\Delta z$ then needs to be projected back to $\Delta x$. Both the solution of the linear system and the projection are conducted twice for Mehrotra's predictor-corrector algorithm.
	
	If there are only bound constraints, forming the sparse matrix product $A_{i}^TW_{i}^{-1}\Lambda_{i}A_{i}$ for the classical normal equations requires only negligible copying operations. The cost of the Cholesky decomposition of $\Phi$ can be neglected in case of a diagonal $H$ and a low rank matrix $A_{i}^TW_{i}^{-1}\Lambda_{i}A_{i}$ allowing for cheap decomposition updates. 
	In case of a non-diagonal $H$ each Newton iteration requires a further Cholesky decomposition of the matrix $\Phi$ requiring $O(T(n_x + n_u)^3/3)$ operations.
	The Cholesky decompositions of $A_e\Phi^{-1}A_e^T$ and $N^T\Phi N$ for the classical and projected normal equations, respectively, are of the same cost. Solving the linear system for the primal and dual is slightly more expensive for the projected normal equations due to the additional virtual control variables.
	
	Depending on the number of inequality constraints $m_i$ and controls $n_u$, the operational counts may shift in favor of FORCESPRO (high $m_i$ and low $n_u$) or NIPM-MPC (low $m_i$ and high $n_u$) but overall are very similar for the case of diagonal $Q$ and $R$ matrices. In~\cite{zeilinger2012} it is not further specified under which circumstances $A_i^TW_i^{-1}\Lambda_iA_i$ can be considered low-rank. In case of bound constraints both on the states and the controls, as it is common for many control applications, the matrix $A_i^TW_i^{-1}\Lambda_iA_i$ is actually of full rank $n$ (assuming that $\Lambda_i\neq0$) and supposedly inapplicable for the update scheme. In this case NIPM-MPC would be in advantage since only a single Cholesky decomposition instead of two is required per Newton iteration. 
	
	\section{algorithm}
	\algnewcommand{\IIf}[1]{\State\algorithmicif\ #1\ \algorithmicthen}
	\algnewcommand{\EndIIf}{\unskip\ \algorithmicend\ \algorithmicif}
	
	\begin{algorithm}[t!]
		\caption{Primal-dual NIPM-MPC}\label{alg:ipmHLSP}
		\begin{algorithmic}[1]
			\Procedure{NIPM-MPC}{$A_iN$, $N^THN$, $QR(A_{e})$}
			\State Solve $A_ex = b_e$ for the initial feasible $x$ with the sparse QR decomposition $QR(A_e)$
			\For{$i < i_{max}$}
			\State Calculate $N^Tr_1$ with predictor $F$~\eqref{eq:fpredictor}
			\IIf{$\Vert K \Vert < \epsilon$ or $\Vert N^Tr_1 \Vert < \epsilon$} Break 	\EndIIf
			\State Calculate $(A_iN)^TW_i^{-1}\Lambda_i(A_iN)$
			\State Cholesky decomposition of $N^T\Phi N$
			\State Solve~\eqref{eq:projNEq} for $\Delta z$
			\State Projection~\eqref{eq:changeOfVarQP}
			\State Line search for dual feasibility~\eqref{eq:feasCond}
			\State Calculate $N^Tr_1$ with corrector $F$~\eqref{eq:fcorrector}
			\State Solve~\eqref{eq:projNEq} for $\Delta z$
			\State Projection~\eqref{eq:changeOfVarQP}
			\State Line search for dual feasibility~\eqref{eq:feasCond} 
			\EndFor
			\EndProcedure
		\end{algorithmic}
	\label{alg:nipmmpc}	
	\end{algorithm}
	
	The overall algorithm of our solver is given in alg.~\ref{alg:nipmmpc}. The projections $A_iN$~\eqref{eq:nain} and $N^THN$~\eqref{eq:nhn} have been computed offline during the solver setup. At the beginning of each Newton's method, the initial feasible equality point $A_ex=b_e$ is obtained. After a convergence test, the composition and Cholesky decomposition of $N^T\Phi N$ is done in one sweep as recommended in~\cite{zeilinger2012} in order to streamline memory access. The decomposition is used to calculate both the predictor and corrector step. The process is repeated until the norm of the KKT vector $K$ or the projected primal residual $N^Tr_1$ is below our desired threshold $\epsilon$. Note that for the calculation of $r_1$ we can use $F$~\eqref{eq:fpredictor} from the predictor step~\eqref{eq:fpredictor} since $\mu\approx 0$ holds at convergence.

	\section{Evaluation}
	\label{sec:evaluation}

	This section serves the purpose of validating our algorithm NIPM-MPC and evaluating it with respect to the three interior point method based solvers FORCES~\cite{zeilinger2012}, FASTMPC~\cite{wangboyd2010} and GUROBI~\cite{gurobi} in terms of computation time per Newton iteration. FORCES and FASTMPC are based on the classical normal equations and make use of the given sparsity of the MPC problem. FORCES is a code generation algorithm tailored to the specific problem constellation and has been shown to be significantly faster than CVXGEN~\cite{CVXGEN} due to a more efficient composition method of the matrix $A_e\Phi^{-1}A_e^T$. Additionally, when the inequality constraint matrix $A_i$ is low-rank a very efficient update scheme for the Cholesky decomposition of $\Phi$ can be applied. For our evaluation we use the commercial software FORCESPRO~\cite{FORCESPro,FORCESNLP} which is derived from the solver FORCES. We furthermore use FASTMPC which is a primal-barrier interior point method with heuristically tuned barrier method only requiring one linear system solve per Newton iteration. GUROBI is a sparse solver but is not specialized for MPC problems. Further technical details are not known about this commercial software. 
	
	The ADMM based solver OSQP~\cite{osqp} is also iterative in nature but typically only requires a single factorization per QP resolution. We therefore reference it only when the overall computation times of computing the control are compared and focus the evaluation of our solver in comparison with other IPM solvers. OSQP is a sparse solver and therefore efficiently solves MPC problems.
	
	All four reference solvers are used at its standard settings. FORCESPRO and OSQP, which are both coded in C and library free, are accessed over their Python interface and the time measurements are done over their incorporated run time measurement tools. The source code of FASTMPC is open-source and based on LAPACK~\cite{lapack}. We compiled its C code with the compiler flag -O3. GUROBI is used over its C++ interface. NIPM-MPC is written in C++ and based on the Eigen linear algebra library~\cite{eigenweb}. All tests are run on an Intel(R) Core(TM) i7-9750H CPU @ 2.60GHz processor with 32 Gb of RAM.
	
	We conduct three simulations with $g=0$ and $\overline{w}=0$. We pick diagonal cost matrices $Q$ and $R$ for the first two simulations. In the first test~\ref{sec:constX} the number of states $n_x$ is fixed but the number of controls $n_u$ and the horizon length $T$ are varied. The second test~\ref{sec:constT} keeps the length of the horizon fixed but the number of variables is increased with a fixed ratio of state to control variables. In the third test both $Q$ and $R$ are dense matrices such that the solvers based on the classical normal equations have to conduct two Cholesky decompositions. Note that while both FORCESPRO and NIPM-MPC consider the diagonal structure of $Q$ and $R$ FASTMPC does not particularly do so and therefore is always conducting two Cholesky decompositions. 
	
	In all simulations and for each solver we sum up the overall solver computation times and divide them by the number of Newton iterations in order to get the computation times per Newton iteration. Naturally, this offsets any overhead from computations that are conducted outside of the Newton method's. Since for our solver obtaining the primal equality point and computing the Lagrange multipliers equates to matrix-vector multiplications we believe that the measurements are representative nonetheless and can be confirmed with a look at the overall computation times. 
	
		\begin{figure}[htp!]
		\includegraphics[width=\columnwidth]{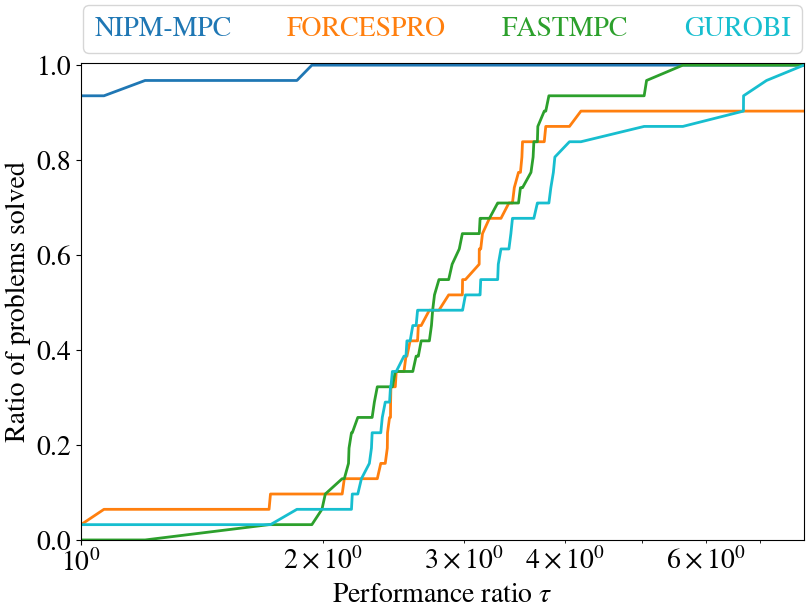}
		\centering
		\caption{Performance profile of the different IPM solvers for the 31 problems solved.}
		\label{fig:perfProfile}
	\end{figure}

	The performance profile of the IPM solvers summarizing the results of the simulations is given in fig.~\ref{fig:perfProfile}. NIPM-MPC is the fastest of the four IPM solvers in 29 of the 31 problems. We choose the required time per Newton iteration as cost for each solver and problem. OSQP is therefore not included.

	\subsection{Constant number of state variables}
	\label{sec:constX}
	
	\begin{figure}[htp!]
		\includegraphics[width=\columnwidth]{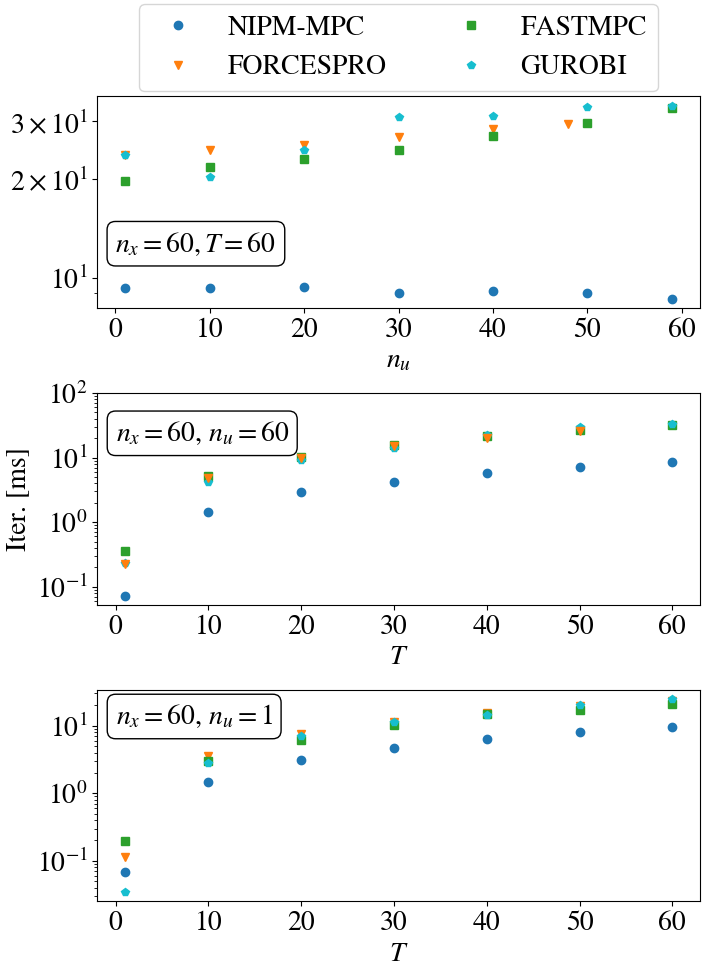}
		\centering
		\caption{Time per Newton iteration for a fixed number of states $n_x$. NIPM-MPC is consistently faster than FORCESPRO and FASTMPC and is  outperformed by GUROBI only for $n_x=60$, $n_u=1$ and $T=1$.}
		\label{fig:randTestX60}
	\end{figure}
	
	First we compare the solvers for a fixed number of states $n_x = 60$. Both matrices $A_{x,e}$ and $B_{u,e}$ are randomly generated dense matrices with spectral radius of one for neutral stability~\cite{wangboyd2010}. We apply bound constraints both on the states $x\in\left[-4,4\right]$ and the controls $u\in\left[-0.5,0.5\right]$. The state vector is initialized uniformly to $0.2$.
	
	The upper graph of fig.~\ref{fig:randTestX60} shows the computation times for an increasing number of controls $n_u$ from $1$ upwards. As expected this does not influence the computation times of our solver NIPM-MPC since the number of original controls $n_u$ and virtual controls $n_{u^*}$ always sums up to $n_{\hat{u}}=n_x$. Consequently, the number of operations is independent of $n_u$ as is outlined in our computational operations table fig.~\ref{fig:NewtonIterOp}. 
	
	For FORCESPRO, FASTMPC and GUROBI the increase of control variables $n_u$ has a direct impact on the computation times as can be seen from fig.~\ref{fig:NewtonIterOp}. The computation times of forming the matrix $A_e\Phi^{-1}A_e^T$ and of computing the primal are linearly and quadratically dependent of $n_u$ as can be recognized from the slight increase in computation times of FORCESPRO and FASTMPC. Timings for FORCESPRO are only given up to $n_u=48$ since for a larger number of controls the code generator fails. For this particular case, NIPM-MPC is around 70\% faster than FORCESPRO. 
	
	This trend continues for a fixed number of controls $n_u = 60$ and $n_u=1$ but increasing length of the receding horizon. As expected, all three solvers show a linear increase in computation times with a higher rate for FORCESPRO, FASTMPC and GUROBI.
	
	\subsection{Fixed ratio of state to control variables and length of receding horizon with diagonal cost matrices}
	\label{sec:constT}
		
	\begin{figure}[htp!]
		\includegraphics[width=\columnwidth]{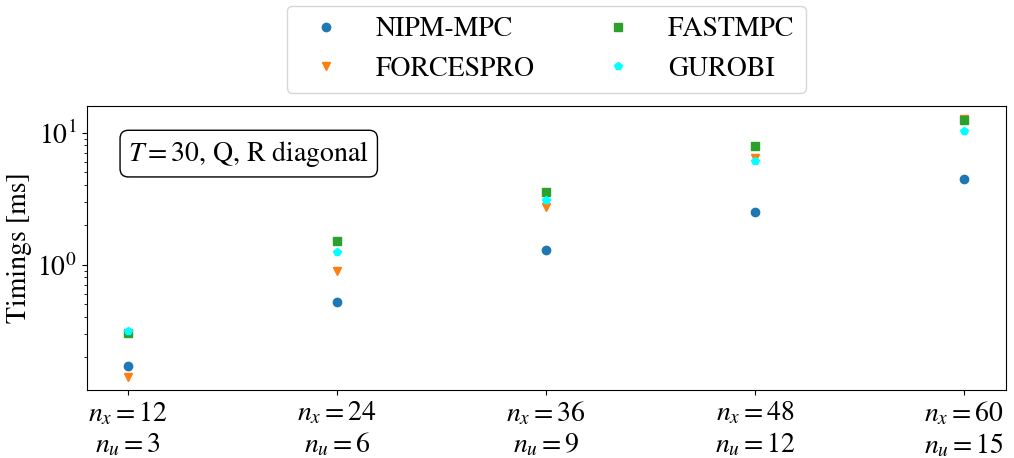}
		\centering
		\caption{Timings for Newton iterations for a fixed ratio of 4:1 between state to control variables. $Q$ and $R$ are diagonal matrices. NIPM-MPC is slower than FORCESPRO for $n_x=12$ but shows a flatter increase in computation times for a growing number of state variables.}
		\label{fig:randTestT60}
	\end{figure}

	In the previous evaluation we considered systems with a high number of states which might be disadvantageous for code generating solvers like FORCESPRO since the code size increases with the problem size~\cite{qpSwift2019}.
	Therefore, we evaluate the solvers now with a lower number of states. The ratio of the number of states with respect to the number of controls is chosen as 4:1. The receding horizon length is $T=30$. Our system and transfer matrices represent the masses-spring dynamic system from~\cite{wangboyd2010}.
	For six masses $M=6$ we have 12 states (six positions and six velocities). Three controls are employed such that our control to states ratio fulfills 4:1. The state vector is initialized uniformly to $1$.
	
	As can be seen in fig.~\ref{fig:randTestT60}, all solvers' computation times increase cubically  with increasing size $n_x$ for the Cholesky decomposition of the matrices $A_e\Phi A_e^T$ (FORCESPRO, FASTMPC) or $\tilde{\Phi}$ (NIPM-MPC). 
	
	For $n_x=12$, the computation times of FORCESPRO are the fastest with $0.14$ms per Newton iteration. NIPM-MPC comes second at $0.17$ms, FASTMPC third at $0.3$ms and GUROBI fourth at $0.31$ms.
	
	\begin{figure}[htp!]
		\includegraphics[width=\columnwidth]{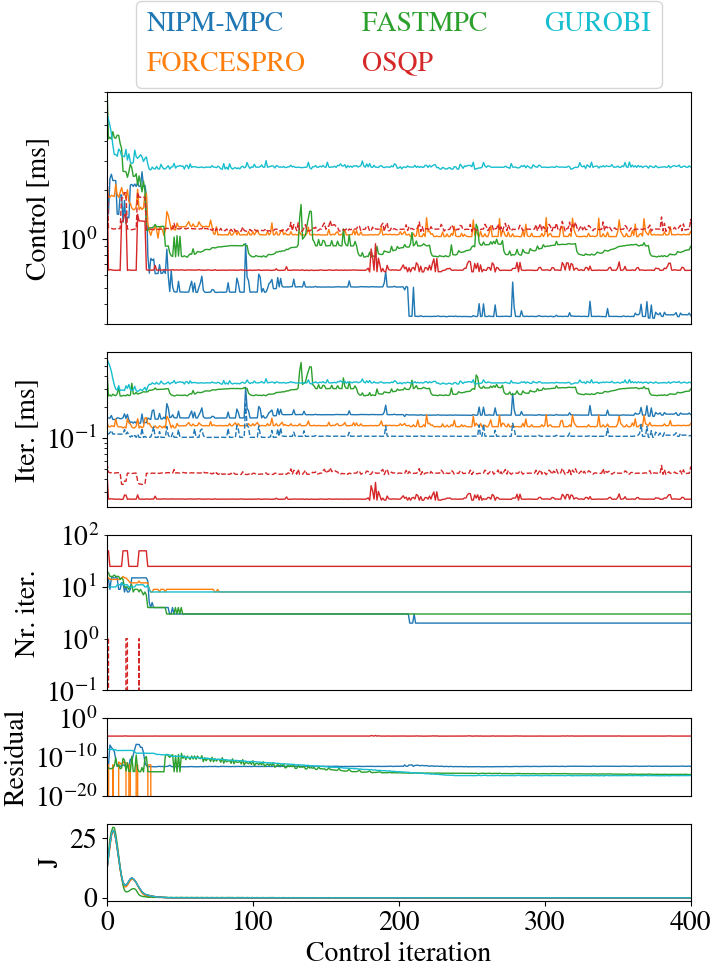}
		\centering
		\caption{Solver and control behavior for mass-spring system with $n_x=12$ and $n_u=3$. Except for the cost $J$ all graphs are in logarithmic scale for better readability.}
		\label{fig:X12U3T30}
	\end{figure}

	In order to further detail these results we plot the solver behaviors for each control step in fig.~\ref{fig:X12U3T30}. The upper graph shows the overall computation times for solving the MPC. The ADMM based solver OSQP generates a moderately accurate solution with a residual of $10^{-4}$ very quickly. It outperforms the IPM based solvers especially in the beginning where all IPM-based solvers require a significant amount of Newton iterations of up to 15 until convergence (see middle graph). However, if a high accuracy solution is desired OSQP requires a significant additional amount of computation since an active-set guess with zero primal and dual residual is generated. The additional computational burden makes it slower than the other specialized IPM based MPC solvers (top graph, dashed red line). Additionally, in control instances where OSQP requires a factorization update (see dashed red line in middle graph) it is approximately as fast as NIPM-MPC. Our solver is thereby reducing the number of Newton iterations necessary in each control instance the quickest to the point that only a single one is necessary and making it faster than the other solvers by a margin.
	
	The second graph from the top indicates the timings of a single solver iteration. FORCESPRO resolves a single Newton iteration the fastest. The dashed blue line shows the contribution of NIPM-MPC's core operations (dashed blue line) of composing $N^T\Phi N$ and its Cholesky decomposition and calculating the projected primal $\Delta z$ twice for Mehrotra's predictor-corrector algorithm. The core operations take less time than a whole FORCESPRO iteration but takes longer considering the rest of the computations of calculating $N^Tr_1$ and conducting the line search.

	\subsection{Fixed ratio of state to control variables and length of receding horizon with dense cost matrices}
	
		\begin{figure}[htp!]
		\includegraphics[width=\columnwidth]{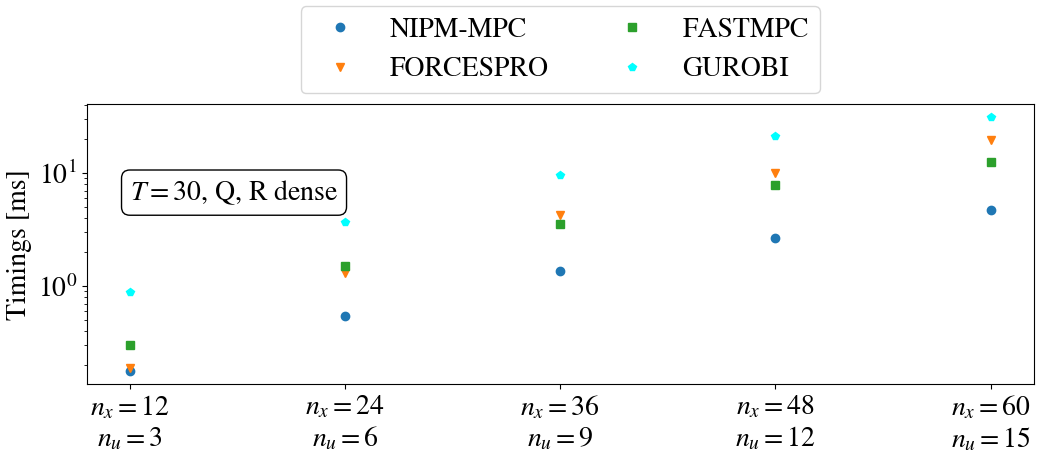}
		\centering
		\caption{Timings for Newton iterations for a fixed ratio of 4:1 between state to control variables. $Q$ and $R$ are dense matrices.}
		\label{fig:randTestT30XU4}
	\end{figure}

	In the last example we want to investigate the behavior of the solvers when $Q$ and $R$ are dense matrices. 
	NIPM-MPC, FORCESPRO and GUROBI handle diagonal cost matrices explicitly and are therefore negatively influenced bu dense ones in terms of computation times.
	FASTMPC always treats $Q$ and $R$ as dense matrices so its computation times do not change. 
	
	As can be seen from fig.~\ref{fig:randTestT30XU4}, all four solver's computational times increase cubically with increasing number of state variables $n_x$. Since now NIPM-MPC, FORCESPRO and GUROBI need to handle dense matrices, their computational times are slightly larger. Explicitly, for $n_x=12$ and $n_u=3$, FORCESPRO takes $0.189$ ms and NIPM-MPC takes $0.177$ ms, making NIPM-MPC slightly faster than FORCESPRO. This is due to the fact that FORCESPRO now needs to conduct an additional Cholesky decomposition of $\Phi$ in each Newton iteration while NIPM-MPC is only influenced by the additional dense matrix vector multiplication $Hy$.
	
	\section{Conclusion}
\label{sec:conclusion}

	In this work we have proposed a new and very efficient IPM to resolve linear receding horizon problems based on the null-space method which requires only a single Cholesky decomposition per Newton iteration instead of two. With our choice of sparse null-space basis and the concept of virtual controls we are able to maintain the block diagonal sparsity of the MPC matrices. Thereby, the cost of the Cholesky decomposition is not affected by increasing numbers of control variables. We showed that our solver is indeed computationally superior with respect to solvers based on classical IPM formulations. Depending on the problem constellation our solver can be up to 70\% faster per Newton iteration. 
	
	In future work we would like to extend our solver to non-linear MPC and explore the possibility of a hierarchical MPC solver which would allow the separation of constraints in a prioritized fashion. An IPM based solver has already been proposed for dense hierarchical least-squares programs~\cite{pfeiffer2021} which requires accurate convergence on each priority level. This may be a distinguishing advantage of the IPM over the ADMM which provides moderately accurate solutions in a very fast fashion.
	

	\bibliographystyle{IEEEtran}
	\bibliography{bib}

\end{document}